\documentclass[12pt]{amsart}

\usepackage{amsmath,amssymb,amsthm,mathtools}
\usepackage[shortlabels]{enumitem}
\usepackage{hyperref}
\usepackage[numbers,sort&compress]{natbib}

\theoremstyle{plain}
\newtheorem{theorem}{Theorem}[section]
\newtheorem{proposition}[theorem]{Proposition}
\newtheorem{lemma}[theorem]{Lemma}
\newtheorem{corollary}[theorem]{Corollary}

\theoremstyle{definition}
\newtheorem{definition}[theorem]{Definition}

\theoremstyle{remark}
\newtheorem{remark}[theorem]{Remark}
\newtheorem{remarks}[theorem]{Remarks}

\numberwithin{equation}{section}
\allowdisplaybreaks
\newcommand{\N}{\mathbb{N}}
\newcommand{\Q}{\mathbb{Q}}
\newcommand{\R}{\mathbb{R}}
\newcommand{\Ninf}{\mathbb{N}_\infty}
\newcommand{\dimH}{\dim_{\mathrm H}}
\newcommand{\dimloc}{\dim_{\mathrm{loc}}}
\newcommand{\Leb}{\lambda}
\newcommand{\SP}{\Sigma_{\mathrm P}}
\newcommand{\SC}{\Sigma_{\mathrm C}}
\newcommand{\fP}{f_{\mathrm P}}
\newcommand{\fC}{f_{\mathrm C}}
\newcommand{\phiP}{\varphi_{\mathrm P}}
\newcommand{\phiC}{\varphi_{\mathrm C}}
\newcommand{\Dif}{\Delta}
\DeclareMathOperator{\diam}{diam}
\DeclareMathOperator{\len}{len}

\begin{document}

\title[Pierce expansions and continued fraction expansions]{A singular homeomorphism between Pierce expansions and continued fraction expansions}

\author{Min Woong Ahn}
\address{Department of Mathematics Education, Dankook University, 152, Jukjeon-ro, Suji-gu, Yongin-si, Gyeonggi-do, 16890, Republic of Korea}
\email{mwahn@dankook.ac.kr}

\date{\today}


\begin{abstract}
In this paper, we study the map $\Theta$ which sends each $x\in[0,1]$ to the number whose continued fraction digits are the first Pierce expansion digit of $x$ and the gaps between two consecutive Pierce expansion digits of $x$. We investigate the regularity of $\Theta$ and $\Theta^{-1}$, including the optimal H\"older exponent of $\Theta^{-1}$, and the local and Hausdorff dimensions of the push-forwards of the Lebesgue measure under $\Theta$ and $\Theta^{-1}$. We also study the differentiability of $\Theta$ at rational and irrational points, and calculate the Hausdorff dimension of the set of points at which $\Theta$ has an infinite derivative.
\end{abstract}

\maketitle
\tableofcontents

\section{Introduction}\label{sec:intro}

Given $x\in[0,1]$, we define
\[
d_1(x)\coloneqq \begin{cases}\lfloor 1/x\rfloor, & \text{if } x\neq0;\\ \infty, & \text{if } x=0,\end{cases}
\qquad\text{and}\qquad
T(x)\coloneqq \begin{cases}1-d_1(x)x, & \text{if } x\neq0;\\ 0, & \text{if } x=0,\end{cases}
\]
respectively, where $\lfloor y\rfloor$ denotes the largest integer not exceeding $y\in\R$. For each $k\in\N$, let $d_k(x)\coloneqq d_1(T^{k-1}(x))$, where $T^{k-1}$ is the $(k-1)$th iterate of $T$ and $T^0$ is the identity map on $[0,1]$. We denote by $\langle d_1(x),d_2(x),\dots\rangle_P$ the Pierce expansion (or alternating Engel expansion) of $x$, i.e.,
\begin{equation}\label{eq:Pierce}
\begin{split}
x=\langle d_1(x),d_2(x),\dots\rangle_P&\coloneqq \sum_{k\in\N}\frac{(-1)^{k+1}}{d_1(x)\cdots d_k(x)}\\
&=\frac{1}{d_1(x)}-\frac{1}{d_1(x)d_2(x)}+\frac{1}{d_1(x)d_2(x)d_3(x)}-\cdots,
\end{split}
\end{equation}
with the conventions $1/\infty\coloneqq 0$ and $c\cdot\infty\coloneqq \infty$ for $c>0$. It is a classical result that the digit sequence $(d_k(x))_{k\in\N}$ is a strictly increasing infinite sequence of positive integers if $x$ is irrational, while it consists of finitely many strictly increasing positive integers followed by $\infty$'s if $x$ is rational. Moreover, via the Pierce expansion, there is a one-to-one correspondence between the irrationals in $[0,1]$ and the strictly increasing infinite sequences of positive integers. We refer the reader to \cite{Pie29,Sha86} and \cite{Ahn23,Ahn24} for the arithmetic and metric properties of Pierce expansions.

In a similar way, for $y\in[0,1]$, we define
\[
a_1(y)\coloneqq \begin{cases}\lfloor 1/y\rfloor, & \text{if } y\neq0;\\ \infty, & \text{if } y=0,\end{cases}
\qquad\text{and}\qquad
S(y)\coloneqq \begin{cases}1/y-a_1(y), & \text{if } y\neq0;\\ 0, & \text{if } y=0,\end{cases}
\]
and $a_k(y)\coloneqq a_1(S^{k-1}(y))$ for each $k\in\N$. Then
\[
y=[a_1(y),a_2(y),\dots]\coloneqq \cfrac{1}{a_1(y)+\cfrac{1}{a_2(y)+\cfrac{1}{\ddots}}}
\]
is the regular continued fraction expansion of $y$. As is well known, the continued fraction expansion gives a one-to-one correspondence between the irrationals in $[0,1]$ and $\N^\N$.

Notice that a strictly increasing infinite sequence $(d_k)_{k\in\N}$ of positive integers is uniquely determined by the first term $d_1$ and the gaps $d_2-d_1,d_3-d_2,\dots$ between consecutive terms, and that the gaps are again positive integers. Hence, it is natural to connect the two expansions in the following simple way. We consider the map $\Theta\colon[0,1]\to[0,1]$ given by
\begin{equation}\label{eq:Theta-intro}
\Theta(x)\coloneqq [d_1(x),\,d_2(x)-d_1(x),\,d_3(x)-d_2(x),\,\dots].
\end{equation}
If $x$ has a finite Pierce expansion of length $n\in\N$, the continued fraction on the right is taken to have exactly $n$ digits. We also put $\Theta(0)\coloneqq0$. We will give a symbolic construction in Section~\ref{sec:aux}; see Definition~\ref{def:FD}. For instance, $\Theta(1/k)=1/k$ for each $k\in\N$, since $1/k=\langle k\rangle_P=[k]$. We remark that metric properties of the gaps between consecutive Pierce digits were already studied by Viader, Bibiloni, and Parad\'is \cite{VBP99}. The representation using the gap sequence $(d_k(x)-d_{k-1}(x))_{k\in\N}$ with $d_0(x)\coloneqq 0$ appears in Pratsiovytyi and Baranovskyi \cite[p.~150]{PB05} (see also \cite[p.~216]{ABPT07}). Albeverio and Torbin \cite{AT15} call it the {\it difference-version} of the Ostrogradsky--Sierpi\'nski--Pierce expansion. In terms of the map $\Theta$, this gap sequence is the continued fraction digit sequence of $\Theta(x)$. More recently, Long, Liang, and Lu \cite{LLL25} studied the Hausdorff dimensions of limsup and liminf sets defined by gap growth. Lu, Long, and Shang \cite{LLS26} established limit theorems for these gaps and determined the Hausdorff dimensions of several exceptional sets defined by their growth.

A brief mention of known singular functions related to $\Theta$ is in order. Recall the Denjoy--Salem formula for Minkowski's question-mark function ${?}\colon[0,1]\to[0,1]$ (see \cite{Den38,Sal43}):
\[
{?}\bigl([a_1,a_2,\dots]\bigr)=\sum_{k\in\N}\frac{(-1)^{k+1}}{2^{a_1+\cdots+a_k-1}}.
\]
In \cite{PVB99}, Parad\'is, Viader, and Bibiloni introduced a continuous, strictly increasing, singular function $F\colon[0,1]\to[0,1]$, which satisfies
\[
F\bigl(\langle d_1,d_2,\dots\rangle_P\bigr)=\sum_{k\in\N}\frac{(-1)^{k+1}}{2^{d_k-1}}.
\]
Notice that the right-hand sides of the above two formulas are nothing but expansions in the alternate dyadic system, in which a number is written as $2^{-m_1}-2^{-m_2}+2^{-m_3}-\cdots$ with $0\le m_1<m_2<\cdots$ (see \cite{PVB01}). Since the partial sums of the gaps in \eqref{eq:Theta-intro} give back the Pierce digits, we have ${?}(\Theta(x))=F(x)$, that is,
\begin{equation}\label{eq:factorization}
\Theta={?}^{-1}\circ F.
\end{equation}
So, $\Theta$ relates Pierce expansions and continued fraction expansions through the alternate dyadic system. We also remark that de Amo, D\'iaz Carrillo, and Fern\'andez-S\'anchez \cite[Theorem~6.3]{dADF22} showed that $F^{-1}=A_\omega$, where $\omega\coloneqq(n/(n+1))_{n\in\N}$ and $A_\omega$ is a singular function obtained by an alternate version of Salem's geometric construction in \cite{Sal43}. Hence, $\Theta^{-1}=F^{-1}\circ{?}=A_\omega\circ{?}$. Notice that the factorization \eqref{eq:factorization} already tells us that $\Theta$ is an increasing homeomorphism, but we will also give a direct proof by using the digit sequences (see Lemma~\ref{lem:homeo}). However, the factorization itself does not give the metric properties and the differentiability of $\Theta$. For example, a composition of two singular functions need not be singular, as ${?}^{-1}\circ{?}$ is the identity map. In this paper, we determine the local and Hausdorff dimensions of the associated measures and the optimal H\"older exponent of $\Theta^{-1}$, and we also study the derivative of $\Theta$ in detail.

Now, we are in a position to state the main results. Throughout the paper, $\Leb$ denotes the Lebesgue measure on $[0,1]$. Put
\[
\mu\coloneqq \Leb\circ\Theta^{-1}\qquad\text{and}\qquad\nu\coloneqq \Leb\circ\Theta,
\]
i.e., $\mu$ and $\nu$ are the push-forwards of $\Leb$ under $\Theta$ and $\Theta^{-1}$, respectively. Then $\mu([0,y])=\Theta^{-1}(y)$ and $\nu([0,x])=\Theta(x)$ for any $x,y\in[0,1]$, i.e., $\Theta^{-1}$ and $\Theta$ are the distribution functions of $\mu$ and $\nu$, respectively. If $m$ is a non-zero finite Borel measure on $[0,1]$, we say that $m$ is {\it carried by} a Borel set $A\subseteq[0,1]$ if $m([0,1]\setminus A)=0$. The {\it Hausdorff dimension of $m$} is defined by
\[
\dimH m\coloneqq \inf\{\dimH A:A\subseteq[0,1] \text{ is a Borel set carrying }m\}
\]
(see \cite[Proposition~10.3]{Fal97}, where this quantity is called the {\em upper Hausdorff dimension of $m$} and is denoted by $\dim_{\mathrm H}^*m$).

The first main result tells us that the map $\Theta$ is indeed a homeomorphism, and that it maps fundamental intervals onto fundamental intervals.

\begin{theorem}\label{thm:homeo}
The map $\Theta\colon[0,1]\to[0,1]$ is an increasing homeomorphism satisfying $\Theta(\Q\cap[0,1])=\Q\cap[0,1]$ and $\Theta(1/k)=1/k$ for all $k\in\N$. Moreover, for each $n\in\N$ and positive integers $d_1<\cdots<d_n$, the map $\Theta$ sends the Pierce fundamental interval with initial digits $(d_1,\dots,d_n)$ onto the continued fraction fundamental interval with initial digits $(d_1,d_2-d_1,\dots,d_n-d_{n-1})$.
\end{theorem}

The second main result is a global quadratic lower bound for the increments of $\Theta$, which is the key ingredient of this paper.

\begin{theorem}\label{thm:quadratic}
For any $x,y\in[0,1]$ with $x\le y$, we have
\[
\Theta(y)-\Theta(x)\ge\frac14(y-x)^2.
\]
Consequently, $|\Theta^{-1}(u)-\Theta^{-1}(v)|\le2|u-v|^{1/2}$ for any $u,v\in[0,1]$, i.e., $\Theta^{-1}$ is H\"older continuous of exponent $1/2$.
\end{theorem}

With the aid of Theorem~\ref{thm:quadratic}, together with some limit theorems in Pierce and continued fraction expansions, we can determine the local dimensions of $\mu$ and $\nu$. For a non-zero finite Borel measure $m$ on $[0,1]$ and $z\in[0,1]$, the {\it local dimension of $m$ at $z$} is defined by
\[
\dimloc m(z)\coloneqq \lim_{r\to0^+}\frac{\log m(B[z,r])}{\log r},
\]
provided that the limit exists, where $B[z,r]\coloneqq [z-r,z+r]\cap[0,1]$ (see \cite[p.~169]{Fal97}). If $\dimloc m(z)=s$ for $m$-almost every $z\in[0,1]$, then we say that $m$ is {\it exact-dimensional with dimension $s$} (see \cite[p.~174]{Fal97}).

\begin{theorem}\label{thm:dimension}
The following hold.
\begin{enumerate}[label = \upshape(\roman*), ref = \roman*, leftmargin=*, widest=iii]
\item \label{itm:dimension-i} $\mu(B[y,r])\le2\sqrt{2r}$ for all $y\in[0,1]$ and $r>0$.
\item \label{itm:dimension-ii} $\mu$ is exact-dimensional with dimension $1/2$, i.e., $\dimloc\mu(y)=1/2$ for $\mu$-almost every $y\in[0,1]$.
\item \label{itm:dimension-iii} $\nu$ is exact-dimensional with dimension $0$, i.e., $\dimloc\nu(x)=0$ for $\nu$-almost every $x\in[0,1]$.
\end{enumerate}
In particular, $\dimH\mu=1/2$ and $\dimH\nu=0$.
\end{theorem}

\begin{corollary}\label{cor:consequences}
The following hold.
\begin{enumerate}[label = \upshape(\roman*), ref = \roman*, leftmargin=*, widest=iv]
\item \label{itm:consequences-i} Both $\mu$ and $\nu$ are singular with respect to $\Leb$. Consequently, $\Theta$ and $\Theta^{-1}$ are singular functions, i.e., $\Theta'(x)=0$ for $\Leb$-almost every $x$, and $(\Theta^{-1})'(y)=0$ for $\Leb$-almost every $y$.
\item \label{itm:consequences-ii} $\Theta'(x)=+\infty$ for $\nu$-almost every $x\in[0,1]$.
\item \label{itm:consequences-iii} $\dimH\Theta^{-1}(A)\le2\dimH A$ for any $A\subseteq[0,1]$. Moreover, $\Theta^{-1}$ is not H\"older continuous of exponent $\alpha$ for any $\alpha>1/2$.
\item \label{itm:consequences-iv} $\Theta$ is not H\"older continuous of exponent $\alpha$ for any $\alpha>0$.
\end{enumerate}
\end{corollary}

\begin{remarks}\label{rems:dimension}
(1) By the law of large numbers for the gaps of Pierce expansion digits due to Lu, Long, and Shang \cite[Theorem~1.1(i)]{LLS26}, we have $a_n(\Theta(x))=d_n(x)-d_{n-1}(x)\to\infty$ as $n\to\infty$ for $\Leb$-almost every $x$. So, $\mu(\mathcal G)=1$, where
\[
\mathcal G\coloneqq\{y\in[0,1]:a_n(y)\to\infty\text{ as }n\to\infty\}.
\]
By Good \cite{Goo41}, we have $\dimH\mathcal G=1/2$. Hence, $\dimH m\le1/2$ for every non-zero finite Borel measure $m$ on $[0,1]$ carried by $\mathcal G$. In this sense, Theorem~\ref{thm:dimension} tells us that $\mu$ attains the largest possible Hausdorff dimension among these measures.

(2) We can also prove $\dimH\nu=0$ by identifying an explicit Borel set of Hausdorff dimension $0$ carrying $\nu$. See Remark~\ref{rem:nu-carrier}.
\end{remarks}

The last main results are concerned with the differentiability of $\Theta$. The first one gives the one-sided derivatives at rational points explicitly.

\begin{theorem}\label{thm:rational}
Let $r\in(0,1)$ be a rational number with Pierce expansion $r=\langle d_1,\dots,d_n\rangle_P$, where $d_1,\dots,d_n\in\N$. Put $Q_n\coloneqq d_1\cdots d_n$, and let $q_n$ be the denominator of $\Theta(r)=[d_1,d_2-d_1,\dots,d_n-d_{n-1}]$. Then both one-sided derivatives of $\Theta$ at $r$ exist and are finite and positive. More precisely,
\[
\Theta'_-(r)=\frac{Q_n}{q_n^2}\quad\text{and}\quad\Theta'_+(r)=\frac{Q_n(d_n-1)}{q_n^2}\qquad\text{if $n$ is odd},
\]
and
\[
\Theta'_+(r)=\frac{Q_n}{q_n^2}\quad\text{and}\quad\Theta'_-(r)=\frac{Q_n(d_n-1)}{q_n^2}\qquad\text{if $n$ is even}.
\]
Furthermore, $\Theta'_+(0)=\Theta'_-(1)=1$.
\end{theorem}

For instance, $\Theta'_-(1/k)=1/k$ and $\Theta'_+(1/k)=(k-1)/k$ for each integer $k\ge2$.

\begin{theorem}\label{thm:irrational}
Let $x\in(0,1)$ be irrational. If $\Theta'(x)$ exists and is finite, then $\Theta'(x)=0$.
\end{theorem}

\begin{corollary}\label{cor:classification}
We have
\[
\{x\in(0,1):\Theta'(x) \text{ exists and } 0<\Theta'(x)<\infty\}=\{1/2\},
\]
and $\Theta'(1/2)=1/2$.
\end{corollary}

We can also determine the Hausdorff dimension of the set where the derivative is infinite.

\begin{theorem}\label{thm:infinite-dimension}
The set
\[
E_\infty\coloneqq \{x\in(0,1):\Theta'(x)=+\infty\}
\]
has Hausdorff dimension $1/2$.
\end{theorem}

\begin{remark}\label{rem:minkowski}
It is known that ${?}'(x)$, whenever it exists, is either $0$ or $+\infty$ (see \cite{PVB01}). Unlike Minkowski's question-mark function, our map $\Theta$ has a finite positive derivative at exactly one point, namely $x=1/2$. At every other rational point in $(0,1)$, the map $\Theta$ has two different finite positive one-sided derivatives, and so $\Theta$ is not differentiable on a dense subset of $[0,1]$. Moreover, by Corollary~\ref{cor:consequences}(\ref{itm:consequences-i}) and (\ref{itm:consequences-ii}), we have $\Theta'=0$ $\Leb$-almost everywhere and $\Theta'=+\infty$ $\nu$-almost everywhere. We refer the reader to \cite{FVPD12,FVPD14} for continuous, strictly increasing, singular functions having a non-zero finite derivative on an uncountable set and on a dense set, respectively.
\end{remark}

This paper is organized as follows. In Section~\ref{sec:prelim}, we introduce some basic notions and facts about Pierce expansions and continued fraction expansions, and some useful tools for the dimension of measures. In Section~\ref{sec:aux}, we give a symbolic construction of the map $\Theta$, together with a family of auxiliary maps $F_D$, $D\in\N_0$, and establish some auxiliary results. In Section~\ref{sec:proofs}, we prove the main results mentioned in this section. In Section~\ref{sec:remarks}, we conclude with a Lochs-type comparison of the two expansions.

Throughout the paper, we denote by $\N$ the set of positive integers, $\N_0$ the set of non-negative integers, and $\Ninf\coloneqq \N\cup\{\infty\}$ the set of extended positive integers. Following the convention, we define $1/\infty\coloneqq 0$ and $\infty\pm c\coloneqq \infty$ for any $c\in\R$. We order $\Ninf$ so that $\infty$ is larger than every positive integer. For an interval $I\subseteq\R$, we write $\len I$ for its length. For a subset $A$ of $[0,1]$, we denote by $\diam A$ the diameter and $\dimH A$ the Hausdorff dimension. The logarithm without base, denoted $\log$, will always mean the natural logarithm.

\section{Preliminaries}\label{sec:prelim}

In this section, we introduce some basic notions and facts about Pierce expansions and continued fraction expansions, and some tools for the dimension of measures, which will be used in the subsequent sections.

\subsection{Pierce and continued fraction expansions}

\subsubsection{Pierce expansions}
Observe that by definition, for each $j\in\N$, we have $d_1(x)=j$ if and only if $x\in(1/(j+1),1/j]$, and in this case, $T(x)\in[0,1/(j+1))$.

We shall consider a symbolic space closely related to the Pierce expansion digit sequences, which was introduced in \cite{Ahn23}. Let $\Sigma_{\mathrm P,0}\coloneqq \{(\infty,\infty,\dots)\}$, and for each $n\in\N$, let
\[
\Sigma_{\mathrm P,n}\coloneqq \{(\sigma_k)_{k\in\N}\in\N^{\{1,\dots,n\}}\times\{\infty\}^{\N\setminus\{1,\dots,n\}}:\sigma_1<\sigma_2<\cdots<\sigma_n\}.
\]
For brevity, we write $(\sigma_1,\dots,\sigma_n)$ for $(\sigma_1,\dots,\sigma_n,\infty,\infty,\dots)\in\Sigma_{\mathrm P,n}$ whenever no confusion could arise. For the set of all strictly increasing infinite sequences of positive integers, put
\[
\Sigma_{\mathrm P,\infty}\coloneqq \{(\sigma_k)_{k\in\N}\in\N^\N:\sigma_k<\sigma_{k+1}\text{ for all }k\in\N\}.
\]
Finally, let
\[
\SP\coloneqq \bigcup_{n\in\N_0}\Sigma_{\mathrm P,n}\cup\Sigma_{\mathrm P,\infty}.
\]
Each element of $\SP$ is said to be a {\it Pierce sequence}, and an element of $\Sigma_{\mathrm P,n}$, $n\in\N_0$, is said to be of {\it length $n$}. We say that $\sigma\coloneqq (\sigma_k)_{k\in\N}\in\SP$ is {\it realizable} if there exists $x\in[0,1]$ such that $d_k(x)=\sigma_k$ for all $k\in\N$, and we denote by $\SP^{\mathrm{re}}$ the collection of all realizable Pierce sequences.

We denote by $\fP\colon[0,1]\to\SP$ the map sending a number to its sequence of Pierce expansion digits, i.e., $\fP(x)\coloneqq (d_k(x))_{k\in\N}$ for each $x\in[0,1]$. Conversely, define a map $\phiP\colon\SP\to[0,1]$ by
\[
\phiP(\sigma)\coloneqq \sum_{k\in\N}\frac{(-1)^{k+1}}{\sigma_1\cdots\sigma_k}
\]
for each $\sigma\coloneqq (\sigma_k)_{k\in\N}\in\SP$, which is modelled on \eqref{eq:Pierce}. For $\sigma\in\Sigma_{\mathrm P,n}$ with $n\in\N$, we also write $Q(\sigma)\coloneqq \sigma_1\cdots\sigma_n$.

The following facts are classical (see \cite[pp.~23--24]{Sha86}), and they are stated in our notation in \cite[Section~2]{Ahn23}.

\begin{proposition}[See {\cite[Propositions~2.1--2.4]{Ahn23}}]\label{prop:Pierce-basic}
The following hold.
\begin{enumerate}[label = \upshape(\roman*), ref = \roman*, leftmargin=*, widest=iii]
\item \label{itm:Pierce-basic-i} For each $x\in[0,1]$, we have $\fP(x)\in\SP^{\mathrm{re}}$ and $\phiP(\fP(x))=x$. Moreover, $\fP(x)\in\Sigma_{\mathrm P,\infty}$ if and only if $x$ is irrational.
\item \label{itm:Pierce-basic-ii} A Pierce sequence $\sigma$ is not realizable if and only if $\sigma\in\Sigma_{\mathrm P,n}$ for some $n\ge2$ with $\sigma_n=\sigma_{n-1}+1$. In addition, $\fP\colon[0,1]\to\SP^{\mathrm{re}}$ is a bijection with the inverse $\phiP|_{\SP^{\mathrm{re}}}$.
\item \label{itm:Pierce-basic-iii} If $d_n(x)<\infty$ for some $n\in\N$, then
\begin{equation}\label{eq:Pierce-tail}
x=\phiP(d_1(x),\dots,d_n(x))+\frac{(-1)^nT^n(x)}{d_1(x)\cdots d_n(x)}\quad\text{with}\quad 0\le T^n(x)<\frac{1}{d_n(x)+1},
\end{equation}
and $\fP(T^n(x))=(d_{n+k}(x))_{k\in\N}$.
\end{enumerate}
\end{proposition}

Let $n\in\N$ and $\sigma\in\Sigma_{\mathrm P,n}\cap\SP^{\mathrm{re}}$ with $\sigma_n\ge2$. Then $\sigma_n-1>\sigma_{n-1}$, where $\sigma_0\coloneqq 0$, and so $(\sigma_1,\dots,\sigma_{n-1},\sigma_n-1,\sigma_n)\in\Sigma_{\mathrm P,n+1}$. Since $1/(\sigma_n-1)-1/((\sigma_n-1)\sigma_n)=1/\sigma_n$, it is clear that
\begin{equation}\label{eq:Pierce-two-reps}
\phiP(\sigma_1,\dots,\sigma_{n-1},\sigma_n-1,\sigma_n)=\phiP(\sigma_1,\dots,\sigma_{n-1},\sigma_n).
\end{equation}
In other words, every rational number in $(0,1)$ has two representations by finite Pierce sequences, one of which is not realizable (see also \cite[Proposition~2.3]{Ahn23}).

For each $n\in\N$ and $\sigma\in\Sigma_{\mathrm P,n}$, the {\it fundamental interval} associated with $\sigma$ is defined by
\[
I_\sigma\coloneqq \{x\in[0,1]:d_k(x)=\sigma_k\text{ for all }1\le k\le n\}.
\]
The following proposition describes $I_\sigma$ precisely.

\begin{proposition}[See {\cite[Proposition~3.5]{Ahn23}}]\label{prop:Pierce-cylinder}
Let $n\in\N$ and $\sigma\in\Sigma_{\mathrm P,n}$. Put $r\coloneqq\phiP(\sigma)$ and $Q\coloneqq Q(\sigma)$. Then
\[
I_\sigma=\left\{r+\frac{(-1)^nt}{Q}:0\le t<\frac{1}{\sigma_n+1},\ t\neq0\text{ if }\sigma\notin\SP^{\mathrm{re}}\right\},
\]
where the parameter corresponding to each $x\in I_\sigma$ is $t=T^n(x)$. In particular, $I_\sigma$ is an interval with $\len I_\sigma=1/(Q\cdot(\sigma_n+1))$.
\end{proposition}

\begin{proof}
The interval description and its length follow from \cite[Proposition~3.5]{Ahn23}. The identification $t=T^n(x)$ follows from \eqref{eq:Pierce-tail}.
\end{proof}

\subsubsection{Continued fraction expansions}
As in the case of Pierce expansions, for each $j\in\N$, we have $a_1(y)=j$ if and only if $y\in(1/(j+1),1/j]$. Note also that $1=[1]$ by definition. We consider a symbolic space for continued fractions, which is similar to the one in \cite{Ahn25}. Let $\Sigma_{\mathrm C,0}\coloneqq \{(\infty,\infty,\dots)\}$, and for each $n\in\N$, let $\Sigma_{\mathrm C,n}\coloneqq \N^{\{1,\dots,n\}}\times\{\infty\}^{\N\setminus\{1,\dots,n\}}$. Put
\[
\SC\coloneqq \bigcup_{n\in\N_0}\Sigma_{\mathrm C,n}\cup\N^\N\qquad\text{and}\qquad\SC^{\mathrm{re}}\coloneqq \SC\setminus\{\tau\in\Sigma_{\mathrm C,n}:n\ge2\text{ and }\tau_n=1\}.
\]
We write $\fC(y)\coloneqq (a_k(y))_{k\in\N}$ for $y\in[0,1]$, and $\phiC(\tau)\coloneqq [\tau_1,\tau_2,\dots]$ for $\tau\coloneqq (\tau_k)_{k\in\N}\in\SC$, with the convention $1/\infty=0$. For $\tau\in\Sigma_{\mathrm C,n}$ with $n\in\N$, the {\it fundamental interval} associated with $\tau$ is defined by
\[
J_\tau\coloneqq \{y\in[0,1]:a_k(y)=\tau_k\text{ for all }1\le k\le n\}.
\]
Recall that $[\tau_1,\dots,\tau_n]=[\tau_1,\dots,\tau_n-1,1]$ whenever $\tau_n\ge2$.

The following facts are standard; see, e.g., \cite{IK02,Khi64}.

\begin{proposition}\label{prop:CF-basic}
The following hold.
\begin{enumerate}[label = \upshape(\roman*), ref = \roman*, leftmargin=*, widest=iii]
\item \label{itm:CF-basic-i} For each $y\in[0,1]$, we have $\fC(y)\in\SC^{\mathrm{re}}$ and $\phiC(\fC(y))=y$. Moreover, $\fC\colon[0,1]\to\SC^{\mathrm{re}}$ is a bijection with the inverse $\phiC|_{\SC^{\mathrm{re}}}$, and $\fC(y)\in\N^\N$ if and only if $y$ is irrational.
\item \label{itm:CF-basic-ii} Let $n\in\N$ and $\tau\in\Sigma_{\mathrm C,n}$, and let $p_k/q_k\coloneqq [\tau_1,\dots,\tau_k]$, $0\le k\le n$, be its convergents, where $p_{-1}\coloneqq 1$, $q_{-1}\coloneqq 0$, $p_0\coloneqq 0$, $q_0\coloneqq 1$, $p_k\coloneqq \tau_kp_{k-1}+p_{k-2}$, and $q_k\coloneqq \tau_kq_{k-1}+q_{k-2}$. Then $J_\tau$ consists of the points
\[
y=\frac{p_n+p_{n-1}u}{q_n+q_{n-1}u},\quad\text{so that}\quad y-\frac{p_n}{q_n}=\frac{(-1)^nu}{q_n(q_n+q_{n-1}u)},
\]
where $u\coloneqq S^n(y)\in[0,1)$, and $u\neq0$ if $\tau\not\in\SC^{\mathrm{re}}$. In particular, $\len J_\tau=1/(q_n(q_n+q_{n-1}))$.
\item \label{itm:CF-basic-iii} With the notation of (\ref{itm:CF-basic-ii}), we have $\prod_{k=1}^n\tau_k\le q_n\le\prod_{k=1}^n(\tau_k+1)$.
\end{enumerate}
\end{proposition}

\subsubsection{An ordering lemma}
For two distinct sequences $\sigma,\tau\in\Ninf^\N$, let $k(\sigma,\tau)$ denote the smallest index $k$ with $\sigma_k\neq\tau_k$. The following lemma says that the order of two numbers is determined by the first different digit, in the same manner for both expansions.

\begin{lemma}\label{lem:order}
Let $x,x'\in[0,1]$ with $x\neq x'$, and put $k\coloneqq k(\fP(x),\fP(x'))$. Then $x<x'$ if and only if either $k$ is odd and $d_k(x)>d_k(x')$, or $k$ is even and $d_k(x)<d_k(x')$. The same statement holds for continued fraction expansions, i.e., with $\fC$ and $a_k$ in place of $\fP$ and $d_k$.
\end{lemma}

\begin{proof}
Since $\fP$ is injective by Proposition~\ref{prop:Pierce-basic}(\ref{itm:Pierce-basic-ii}), $k$ is well defined. We claim that
\begin{equation}\label{eq:order-diff}
x-x'=\frac{(-1)^{k-1}}{d_1(x)\cdots d_{k-1}(x)}\bigl(T^{k-1}(x)-T^{k-1}(x')\bigr),
\end{equation}
where the empty product is $1$ and $T^0$ is the identity map. If $k=1$, then \eqref{eq:order-diff} is trivial. If $k\ge2$, then $d_{k-1}(x)=d_{k-1}(x')$ is finite, since otherwise both digit sequences would be $\infty$ from the $(k-1)$th term on, contradicting $d_k(x)\neq d_k(x')$. So, \eqref{eq:order-diff} follows from \eqref{eq:Pierce-tail} with $n=k-1$, applied to $x$ and $x'$.
Put $t\coloneqq T^{k-1}(x)$ and $t'\coloneqq T^{k-1}(x')$. Then $d_1(t)=d_k(x)\neq d_k(x')=d_1(t')$. Notice that the sets $\{s\in[0,1]:d_1(s)=j\}=(1/(j+1),1/j]$, $j\in\N$, and $\{s\in[0,1]:d_1(s)=\infty\}=\{0\}$ are pairwise disjoint and arranged in the decreasing order of $j\in\Ninf$. So, $d_k(x)>d_k(x')$ implies $t<t'$, and the assertion follows.

For continued fraction expansions, we have $y=(p_{k-1}+p_{k-2}u)/(q_{k-1}+q_{k-2}u)$ with $u=S^{k-1}(y)$ by Proposition~\ref{prop:CF-basic}(\ref{itm:CF-basic-ii}) if $k\ge2$; this also holds trivially if $k=1$, since then $y=u$ by $p_0=0$, $q_0=1$, $p_{-1}=1$, and $q_{-1}=0$. The map $u\mapsto(p_{k-1}+p_{k-2}u)/(q_{k-1}+q_{k-2}u)$ is strictly monotone, and it is increasing if and only if $p_{k-2}q_{k-1}-p_{k-1}q_{k-2}=(-1)^{k-1}>0$. Since the first-digit sets are the same as above, the same argument applies.
\end{proof}

\subsection{Dimension results and limit theorems}
We will use the following two well-known facts.

\begin{proposition}[See {\cite[Proposition~10.3]{Fal97}}]\label{prop:local-dim}
Let $m$ be a non-zero finite Borel measure on $[0,1]$. If $\dimloc m(z)=s$ for $m$-almost every $z\in[0,1]$, then $\dimH m=s$.
\end{proposition}

\begin{proposition}[See {\cite[Proposition~3.3(a)]{Fal14}}]\label{prop:holder-dim}
Let $A\subseteq[0,1]$ and $g\colon A\to\R$. Suppose that there exist constants $c>0$ and $\alpha\in(0,1]$ such that $|g(x)-g(y)|\le c|x-y|^\alpha$ for all $x,y\in A$. Then $\dimH g(A)\le\alpha^{-1}\dimH A$.
\end{proposition}

We will also use the following three limit theorems.
\begin{enumerate}[label = \upshape(L\arabic*), ref = L\arabic*, leftmargin=*, widest=L3]
\item \label{lim:L1} (Shallit \cite[Theorem~16]{Sha86}) We have $(\log d_n(x))/n\to1$ as $n\to\infty$ for $\Leb$-almost every $x\in[0,1]$.
\item \label{lim:L2} (Lu, Long, and Shang \cite[Theorem~1.1(i)]{LLS26}) For $\Leb$-almost every $x\in[0,1]$, we have $(\log(d_n(x)-d_{n-1}(x)))/n\to1$ as $n\to\infty$, where $d_0(x)\coloneqq 0$.
\item \label{lim:L3} (L\'evy \cite{Lev36}) For $\Leb$-almost every $y\in[0,1]$, we have $(\log q_n(y))/n\to\pi^2/(12\log2)$ as $n\to\infty$, where $q_n(y)$ is the denominator of the $n$th convergent $[a_1(y),\dots,a_n(y)]$.
\end{enumerate}

\section{The map \texorpdfstring{$\Theta$}{Theta} and auxiliary results}\label{sec:aux}

In this section, we give a symbolic construction of the map $\Theta$ given in \eqref{eq:Theta-intro}. In fact, it is convenient to consider a family of auxiliary maps $F_D$, $D\in\N_0$, containing $\Theta=F_0$, since the tail of a Pierce sequence again gives such a map with larger $D$.

For each $D\in\N_0$, put $\Omega_D\coloneqq [0,1/(D+1)]$. Notice that every finite Pierce expansion digit of $t\in\Omega_D$ is larger than $D$, because $d_1(t)\ge D+1$ and the digits are strictly increasing. We also put
\[
\SP^{(D)}\coloneqq \{\sigma\in\SP:\sigma_1>D\},
\]
so that $\fP(\Omega_D)=\SP^{\mathrm{re}}\cap\SP^{(D)}$.

\begin{definition}\label{def:FD}
Let $D\in\N_0$. For each $\sigma\coloneqq(\sigma_k)_{k\in\N}\in\SP^{(D)}$, put $\sigma_0\coloneqq D$ and define
\[
\Dif_D(\sigma)\coloneqq(\tau_k)_{k\in\N}, \qquad \tau_k\coloneqq
\begin{cases}
\sigma_k-\sigma_{k-1}, & \text{if }\sigma_k<\infty,\\
\infty, & \text{if }\sigma_k=\infty.
\end{cases}
\]
Then $\Dif_D(\sigma)\in\SC$. Define $F_D\colon\Omega_D\to[0,1]$ by
\[
F_D(t)\coloneqq\phiC\bigl(\Dif_D(\fP(t))\bigr)
\]
for each $t\in\Omega_D$. In particular, we put $\Theta\coloneqq F_0$.
\end{definition}

In other words, if $t=\langle d_1,d_2,\dots\rangle_P\in\Omega_D$, then $F_D(t)=[d_1-D,\,d_2-d_1,\,d_3-d_2,\dots]$. Clearly, $F_D(0)=0$, and $F_D(1/k)=1/(k-D)$ for each integer $k\ge D+1$.

\begin{lemma}\label{lem:Dif}
Let $D\in\N_0$. Then the following hold.
\begin{enumerate}[label = \upshape(\roman*), ref = \roman*, leftmargin=*, widest=ii]
\item \label{itm:Dif-i} The map $\Dif_D$ is a bijection from $\SP^{\mathrm{re}}\cap\SP^{(D)}$ onto $\SC^{\mathrm{re}}$. Consequently, $\fC(F_D(t))=\Dif_D(\fP(t))$ for each $t\in\Omega_D$.
\item \label{itm:Dif-ii} For each $n\in\N$ and $\sigma\in\Sigma_{\mathrm P,n}\cap\SP^{(D)}$, realizable or not, we have $F_D(\phiP(\sigma))=\phiC(\Dif_D(\sigma))$.
\end{enumerate}
\end{lemma}

\begin{proof}
(\ref{itm:Dif-i}) It is clear that $\Dif_D$ maps $\Sigma_{\mathrm P,\infty}\cap\SP^{(D)}$ onto $\N^\N$ bijectively, with the inverse $\tau\mapsto(D+\tau_1,\,D+\tau_1+\tau_2,\,\dots)$, and that $\Dif_D$ maps $\Sigma_{\mathrm P,n}\cap\SP^{(D)}$ onto $\Sigma_{\mathrm C,n}$ bijectively for each $n\in\N_0$. Let $n\ge2$ and $\sigma\in\Sigma_{\mathrm P,n}\cap\SP^{(D)}$. By Proposition~\ref{prop:Pierce-basic}(\ref{itm:Pierce-basic-ii}), $\sigma$ is realizable if and only if $\sigma_n-\sigma_{n-1}\ge2$, i.e., if and only if the last finite term of $\Dif_D(\sigma)$ is at least $2$, which means that $\Dif_D(\sigma)\in\SC^{\mathrm{re}}$. Since every sequence of length $0$ or $1$ is realizable in both $\SP$ and $\SC$, this proves the first assertion. The second assertion follows from the first one and Proposition~\ref{prop:CF-basic}(\ref{itm:CF-basic-i}).

(\ref{itm:Dif-ii}) If $\sigma$ is realizable, then $\sigma=\fP(\phiP(\sigma))$ by Proposition~\ref{prop:Pierce-basic}(\ref{itm:Pierce-basic-ii}), and the assertion is nothing but the definition of $F_D$. Assume that $\sigma$ is not realizable. Then $n\ge2$ and $\sigma_n=\sigma_{n-1}+1$ by Proposition~\ref{prop:Pierce-basic}(\ref{itm:Pierce-basic-ii}). Put $\sigma'\coloneqq (\sigma_1,\dots,\sigma_{n-2},\sigma_n)\in\Sigma_{\mathrm P,n-1}$. If $n\ge3$, then $\sigma_n\ge\sigma_{n-2}+2$, and so $\sigma'$ is realizable; if $n=2$, then $\sigma'=(\sigma_2)$ is clearly realizable. By \eqref{eq:Pierce-two-reps} applied to $\sigma'$, we have $\phiP(\sigma)=\phiP(\sigma')$. On the continued fraction side, put $b\coloneqq \sigma_{n-1}-\sigma_{n-2}$, where $\sigma_0\coloneqq D$. Then
\[
\phiC(\Dif_D(\sigma))=[\dots,b,1]=[\dots,b+1]=\phiC(\Dif_D(\sigma')).
\]
Since $\sigma'$ is realizable, it follows that $F_D(\phiP(\sigma))=F_D(\phiP(\sigma'))=\phiC(\Dif_D(\sigma'))=\phiC(\Dif_D(\sigma))$, as desired.
\end{proof}

\begin{lemma}\label{lem:homeo}
Let $D\in\N_0$. Then $F_D$ is an increasing homeomorphism of $\Omega_D$ onto $[0,1]$, and it maps the rationals in $\Omega_D$ onto the rationals in $[0,1]$. Moreover, for each $n\in\N$ and $\sigma\in\Sigma_{\mathrm P,n}\cap\SP^{(D)}$, we have
\begin{equation}\label{eq:cylinder-corr}
F_D(I_\sigma)=J_{\Dif_D(\sigma)}.
\end{equation}
\end{lemma}

\begin{proof}
By Lemma~\ref{lem:Dif}(\ref{itm:Dif-i}) and Propositions~\ref{prop:Pierce-basic}(\ref{itm:Pierce-basic-ii}) and \ref{prop:CF-basic}(\ref{itm:CF-basic-i}), the map $F_D=\phiC\circ\Dif_D\circ\fP$ is a bijection of $\Omega_D$ onto $[0,1]$. Since finite sequences correspond to finite sequences under $\Dif_D$, it follows from Propositions~\ref{prop:Pierce-basic}(\ref{itm:Pierce-basic-i}) and \ref{prop:CF-basic}(\ref{itm:CF-basic-i}) that $F_D$ maps rationals onto rationals.

We show that $F_D$ is increasing. Let $t,t'\in\Omega_D$ with $t<t'$, and write $d\coloneqq \fP(t)$ and $d'\coloneqq \fP(t')$. Put $k\coloneqq k(d,d')$. Then it is clear that $k$ is also the smallest index at which $\Dif_D(d)$ and $\Dif_D(d')$ differ. If both $d_k$ and $d'_k$ are finite, then $\Dif_D(d)_k-\Dif_D(d')_k=d_k-d'_k$. If exactly one of them is $\infty$, say $d_k=\infty$, then $\Dif_D(d)_k=\infty$ while $\Dif_D(d')_k$ is finite. In either case, the order between the $k$th terms is preserved under $\Dif_D$. Since $\fC(F_D(t))=\Dif_D(d)$ and $\fC(F_D(t'))=\Dif_D(d')$ by Lemma~\ref{lem:Dif}(\ref{itm:Dif-i}), we deduce from Lemma~\ref{lem:order}, applied to both expansions, that $F_D(t)<F_D(t')$.

Now, a strictly increasing bijection between two compact intervals is a homeomorphism. Finally, \eqref{eq:cylinder-corr} follows from Lemma~\ref{lem:Dif}(\ref{itm:Dif-i}) and the bijectivity of $F_D$. This completes the proof.
\end{proof}

For each $j\in\N$, let $U_j\coloneqq [1/(j+1),1/j]$. The following recursive formula connects $F_D$ with $F_j$, and it will be the main tool in the proof of Theorem~\ref{thm:quadratic}.

\begin{lemma}\label{lem:recursion}
Let $D\in\N_0$ and $j\in\N$ with $j\ge D+1$, and put $a\coloneqq j-D$. Then, for each $t\in U_j$, we have $1-jt\in\Omega_j$ and
\begin{equation}\label{eq:recursion}
F_D(t)=\frac{1}{a+F_j(1-jt)}.
\end{equation}
In particular, $a_1(F_D(t))=a$ for each $t\in(1/(j+1),1/j]$.
\end{lemma}

\begin{proof}
Let $t\in U_j$. Clearly, $1-jt\in[0,1/(j+1)]=\Omega_j$. Suppose first that $t\in(1/(j+1),1/j]$. Then $d_1(t)=j$, and $T(t)=1-jt$ has the digits $d_2(t),d_3(t),\dots$ by Proposition~\ref{prop:Pierce-basic}(\ref{itm:Pierce-basic-iii}). Hence,
\[
\begin{split}
F_D(t)&=[a,\,d_2(t)-j,\,d_3(t)-d_2(t),\,\dots]\\
&=\frac{1}{a+[d_2(t)-j,\,d_3(t)-d_2(t),\,\dots]}=\frac{1}{a+F_j(T(t))}.
\end{split}
\]
Now suppose $t=1/(j+1)$. Then $1-jt=1/(j+1)$ and $F_j(1/(j+1))=1$, while $F_D(1/(j+1))=1/(j+1-D)=1/(a+1)$. Hence the lemma.
\end{proof}

\begin{lemma}\label{lem:ratio}
The following hold.
\begin{enumerate}[label = \upshape(\roman*), ref = \roman*, leftmargin=*, widest=iii]
\item \label{itm:ratio-i} If $D\in\N$, then $t\le F_D(t)\le(D+1)t$ for all $t\in\Omega_D$.
\item \label{itm:ratio-ii} $\Theta(t)\ge t/2$ for all $t\in[0,1]$.
\item \label{itm:ratio-iii} For each $D\in\N_0$, we have $F_D(t)/t\to1$ as $t\to0^+$.
\end{enumerate}
\end{lemma}

\begin{proof}
Assertions (\ref{itm:ratio-i}) and (\ref{itm:ratio-ii}) are trivial if $t=0$, and (\ref{itm:ratio-iii}) concerns only $t>0$. So, let $t\neq0$. Then $t\in U_j$ for some $j\ge D+1$; put $a\coloneqq j-D$. Since $F_j$ takes values in $[0,1]$, it follows from \eqref{eq:recursion} that
\begin{equation}\label{eq:FD-range}
\frac{1}{a+1}\le F_D(t)\le\frac1a.
\end{equation}

(\ref{itm:ratio-i}) Let $D\in\N$. Since $j\ge a+1$, we have $F_D(t)\ge1/(a+1)\ge1/j\ge t$ by \eqref{eq:FD-range}. For the upper bound, suppose first that $a\ge2$. Then $(D+1)t\ge(D+1)/(a+D+1)\ge1/a$, where the last inequality holds because $a(D+1)-(a+D+1)=(a-1)D-1\ge0$. Hence, $F_D(t)\le1/a\le(D+1)t$. Suppose next that $a=1$, i.e., $j=D+1$. Put $v\coloneqq 1-(D+1)t\in\Omega_{D+1}$, so that $(D+1)t=1-v$ and $F_D(t)=1/(1+F_{D+1}(v))$ by \eqref{eq:recursion}. Then it suffices to show that $F_{D+1}(v)\ge v/(1-v)$. This is trivial if $v=0$. If $v\neq0$, then $v\in U_{j'}$ for some $j'\ge D+2$, and so, with $a'\coloneqq j'-D-1$, we obtain
\[
F_{D+1}(v)\ge\frac{1}{a'+1}=\frac{1}{j'-D}\ge\frac{1}{j'-1}=\frac{1/j'}{1-1/j'}\ge\frac{v}{1-v}
\]
by \eqref{eq:FD-range}, where we used $D\ge1$, $v\le1/j'$, and the fact that $x\mapsto x/(1-x)$ is increasing on $[0,1)$.

(\ref{itm:ratio-ii}) In this case, $D=0$ and $a=j$. So, by \eqref{eq:FD-range}, $\Theta(t)\ge1/(j+1)\ge1/(2j)\ge t/2$.

(\ref{itm:ratio-iii}) Since $t\in U_j$ and $F_D(t)\in[1/(j-D+1),1/(j-D)]$ by \eqref{eq:FD-range}, we have
\[
\frac{j}{j-D+1}\le\frac{F_D(t)}{t}\le\frac{j+1}{j-D}.
\]
Notice that $j\to\infty$ as $t\to0^+$, and both bounds tend to $1$ as $j\to\infty$. This completes the proof.
\end{proof}

\begin{lemma}\label{lem:branch}
Let $D\in\N_0$ and $j\in\N$ with $j\ge D+1$, and put $a\coloneqq j-D$. Let $t\in U_j$, and put $h\coloneqq 1/j-t$ and $h'\coloneqq t-1/(j+1)$. Then
\begin{align}
F_D(1/j)-F_D(t)&\ge\frac{j}{a(a+1)}h, \label{eq:branch-right}\\
F_D(t)-F_D(1/(j+1))&\ge\frac{j(j+1)}{(a+1)^2}h'. \label{eq:branch-left}
\end{align}
Moreover, the following hold.
\begin{enumerate}[label = \upshape(\roman*), ref = \roman*, leftmargin=*, widest=ii]
\item \label{itm:branch-i} If $D\in\N$, then the right-hand sides of \eqref{eq:branch-right} and \eqref{eq:branch-left} are at least $6h^2$ and $9h'^2$, respectively.
\item \label{itm:branch-ii} If $D=0$, then the right-hand sides of \eqref{eq:branch-right} and \eqref{eq:branch-left} are at least $jh^2$ and $j^2h'^2$, respectively.
\end{enumerate}
\end{lemma}

\begin{proof}
Put $v\coloneqq 1-jt\in\Omega_j$ and $w\coloneqq F_j(v)\in[0,1]$, so that $F_D(t)=1/(a+w)$ by \eqref{eq:recursion}.

To see \eqref{eq:branch-right}, note that $F_D(1/j)=1/a$ and $v=jh$. Since $j\ge1$, Lemma~\ref{lem:ratio}(\ref{itm:ratio-i}) tells us that $w\ge v$. Hence,
\[
F_D(1/j)-F_D(t)=\frac1a-\frac{1}{a+w}=\frac{w}{a(a+w)}\ge\frac{w}{a(a+1)}\ge\frac{jh}{a(a+1)}.
\]
To see \eqref{eq:branch-left}, note that $F_D(1/(j+1))=1/(a+1)$ and $v=1/(j+1)-jh'$. By Lemma~\ref{lem:ratio}(\ref{itm:ratio-i}) again, $w\le(j+1)v=1-j(j+1)h'$. Hence,
\[
F_D(t)-F_D(1/(j+1))=\frac{1}{a+w}-\frac{1}{a+1}=\frac{1-w}{(a+w)(a+1)}\ge\frac{1-w}{(a+1)^2}\ge\frac{j(j+1)h'}{(a+1)^2}.
\]

For the second assertion, observe that $h,h'\le\len U_j=1/(j(j+1))$, i.e., $h\ge j(j+1)h^2$ and $h'\ge j(j+1)h'^2$. Then the right-hand sides of \eqref{eq:branch-right} and \eqref{eq:branch-left} are at least
\[
\frac{j^2(j+1)}{a(a+1)}h^2\qquad\text{and}\qquad\frac{j^2(j+1)^2}{(a+1)^2}h'^2,
\]
respectively. If $D\in\N$, then $j\ge a+1$, and so
\[
\frac{j^2(j+1)}{a(a+1)}\ge\frac{(a+1)(a+2)}{a}=a+3+\frac2a\ge6\qquad\text{and}\qquad\frac{j^2(j+1)^2}{(a+1)^2}\ge(a+2)^2\ge9,
\]
which proves (\ref{itm:branch-i}). If $D=0$, then $a=j$, and the two coefficients become $j$ and $j^2$, which proves (\ref{itm:branch-ii}).
\end{proof}

\section{Proofs of the main results}\label{sec:proofs}

In this section, we prove the main results of this paper mentioned in Section~\ref{sec:intro}.

\begin{proof}[Proof of Theorem~\ref{thm:homeo}]
By Definition~\ref{def:FD}, we have $\Theta=F_0$. Hence, the theorem follows from Lemma~\ref{lem:homeo} with $D=0$, together with $\Theta(1/k)=1/k$ for each $k\in\N$.
\end{proof}

\subsection{The quadratic lower bound}
We first prove the following proposition, which is a version of Theorem~\ref{thm:quadratic} for $F_D$ with $D\in\N$, and which has the better constant $1$.

\begin{proposition}\label{prop:FD-quadratic}
Let $D\in\N$. Then, for any $x,y\in\Omega_D$ with $x\le y$, we have $F_D(y)-F_D(x)\ge(y-x)^2$.
\end{proposition}

Before the proof, we fix some notation which will also be used in the proof of Theorem~\ref{thm:quadratic}. Let $D\in\N_0$ and $x,y\in\Omega_D$ with $x<y$. Since $(0,1/(D+1)]$ is covered by the intervals $U_j$, $j\ge D+1$, one of the following three cases occurs.

{\sc Case 1.} $x=0$.

{\sc Case 2.} $x>0$, and there is no $j$ such that both $x$ and $y$ lie in $U_j$. In this case, take $j_1,j_2\ge D+1$ with $x\in U_{j_1}$ and $y\in U_{j_2}$. Then $j_1\neq j_2$. Moreover, $j_1<j_2$ is impossible, since otherwise $x\ge1/(j_1+1)\ge1/j_2\ge y$. Hence, $j_1>j_2$, and so $x\le1/j_1\le1/(j_2+1)\le y$. Put
\[
h_1\coloneqq \frac{1}{j_1}-x,\qquad h_2\coloneqq \frac{1}{j_2+1}-\frac{1}{j_1},\qquad h_3\coloneqq y-\frac{1}{j_2+1}.
\]
Then $h_1+h_2+h_3=y-x$, and $F_D(y)-F_D(x)$ is the sum of the increments of $F_D$ over the three intervals $[x,1/j_1]$, $[1/j_1,1/(j_2+1)]$, and $[1/(j_2+1),y]$.

{\sc Case 3.} $x,y\in U_j$ for some $j\ge D+1$. In this case, put $a\coloneqq j-D$, $v_x\coloneqq 1-jx$, and $v_y\coloneqq 1-jy$. Then $v_y<v_x$, both $v_x$ and $v_y$ lie in $\Omega_j$, and $v_x-v_y=j(y-x)$. By \eqref{eq:recursion}, we have
\begin{equation}\label{eq:case3}
F_D(y)-F_D(x)=\frac{F_j(v_x)-F_j(v_y)}{(a+F_j(v_x))(a+F_j(v_y))}\ge\frac{F_j(v_x)-F_j(v_y)}{(a+1)^2}.
\end{equation}

\begin{proof}[Proof of Proposition~\ref{prop:FD-quadratic}]
Let $D\in\N$ and $x,y\in\Omega_D$ with $x<y$. We consider the three cases above.

In Case 1, Lemma~\ref{lem:ratio}(\ref{itm:ratio-i}) gives $F_D(y)-F_D(x)=F_D(y)\ge y\ge y^2$.

In Case 2, by Lemma~\ref{lem:branch}(\ref{itm:branch-i}), the increments over the first and the last intervals are at least $6h_1^2$ and $9h_3^2$, respectively. For the middle interval, note that for any integers $k,m$ with $D<k\le m$, we have
\[
\frac{1}{k-D}-\frac{1}{m-D}=\frac{m-k}{(k-D)(m-D)}\ge\frac{m-k}{km}=\frac1k-\frac1m.
\]
Since $F_D(1/k)=1/(k-D)$, applying this with $k=j_2+1$ and $m=j_1$ shows that the increment over the middle interval is at least $h_2$, which is at least $2h_2^2$ because $h_2\le1/2$. Then, by the Cauchy--Schwarz inequality,
\[
(y-x)^2=(h_1+h_2+h_3)^2\le\left(\frac16+\frac12+\frac19\right)\bigl(6h_1^2+2h_2^2+9h_3^2\bigr)=\frac79\bigl(6h_1^2+2h_2^2+9h_3^2\bigr),
\]
and thus $F_D(y)-F_D(x)\ge(9/7)(y-x)^2\ge(y-x)^2$.

For Case 3, we argue by induction on the length of Pierce expansions. For each $N\in\N$, let $\mathcal R_N$ be the set of points in $[0,1]$ whose Pierce expansions have at most $N$ finite digits, and let $P(N)$ be the statement that $F_D(y)-F_D(x)\ge(y-x)^2$ for all $D\in\N$ and all $x,y\in\Omega_D\cap\mathcal R_N$ with $x<y$. Since Cases 1 and 2 are already done for arbitrary points, we only need to consider Case 3 in proving $P(N)$.

We first prove $P(1)$. Note that $\mathcal R_1=\{0\}\cup\{1/k:k\in\N\}$. In Case 3, we must have $x=1/(j+1)$ and $y=1/j$, and so
\[
F_D(y)-F_D(x)=\frac1a-\frac1{a+1}=\frac{1}{a(a+1)}\ge\frac{1}{j^2(j+1)^2}=(y-x)^2,
\]
since $a\le j$.

Now, let $N\ge2$, and assume that $P(N-1)$ holds. Let $x,y\in\mathcal R_N\cap U_j$ with $x<y$. If $x\in(1/(j+1),1/j]$, then $v_x=T(x)$ has at most $N-1$ finite digits by Proposition~\ref{prop:Pierce-basic}(\ref{itm:Pierce-basic-iii}). If $x=1/(j+1)$, then $v_x=1/(j+1)\in\mathcal R_1$. In either case, $v_x\in\mathcal R_{N-1}$, and the same holds for $v_y$. Since $j\ge D+1\ge2$, we can apply $P(N-1)$ to $F_j$, and obtain $F_j(v_x)-F_j(v_y)\ge(v_x-v_y)^2=j^2(y-x)^2$. Hence, by \eqref{eq:case3},
\[
F_D(y)-F_D(x)\ge\frac{j^2}{(a+1)^2}(y-x)^2\ge(y-x)^2,
\]
where the last inequality follows from $j=a+D\ge a+1$. This proves $P(N)$.

By induction, for each $D\in\N$, the desired inequality holds on $\Omega_D\cap\bigcup_{N\in\N}\mathcal R_N=\Q\cap\Omega_D$, which is dense in $\Omega_D$. Since $F_D$ is continuous by Lemma~\ref{lem:homeo}, the inequality holds for all $x,y\in\Omega_D$ with $x\le y$. This completes the proof.
\end{proof}

\begin{proof}[Proof of Theorem~\ref{thm:quadratic}]
Let $x,y\in[0,1]$ with $x<y$. We use the three cases above with $D=0$.

In Case 1, Lemma~\ref{lem:ratio}(\ref{itm:ratio-ii}) gives $\Theta(y)\ge y/2\ge y^2/4$.

In Case 2, we have $j_1\ge2$ and $j_2\ge1$. By Lemma~\ref{lem:branch}(\ref{itm:branch-ii}), the increments of $\Theta$ over the first and the last intervals are at least $j_1h_1^2\ge2h_1^2$ and $j_2^2h_3^2\ge h_3^2$, respectively. The increment over the middle interval is equal to $h_2$, which is at least $h_2^2$, since $\Theta(1/k)=1/k$ for each $k\in\N$. So, by the Cauchy--Schwarz inequality,
\[
(y-x)^2\le\left(\frac12+1+1\right)\bigl(2h_1^2+h_2^2+h_3^2\bigr)\le\frac52\bigl(\Theta(y)-\Theta(x)\bigr),
\]
and thus $\Theta(y)-\Theta(x)\ge(2/5)(y-x)^2\ge(1/4)(y-x)^2$.

In Case 3, we have $a=j\ge1$, and Proposition~\ref{prop:FD-quadratic} applied to $F_j$, together with \eqref{eq:case3}, yields
\[
\Theta(y)-\Theta(x)\ge\frac{j^2}{(j+1)^2}(y-x)^2\ge\frac14(y-x)^2.
\]
This proves the first assertion. For the second one, apply the first one to $\Theta^{-1}(u)$ and $\Theta^{-1}(v)$.
\end{proof}

\subsection{Dimensions of the associated measures}
For irrational $x\in[0,1]$ and $y\coloneqq \Theta(x)$, we write
\[
I_n(x)\coloneqq I_{(d_1(x),\dots,d_n(x))}\qquad\text{and}\qquad J_n(y)\coloneqq J_{(a_1(y),\dots,a_n(y))}
\]
for each $n\in\N$. By Lemma~\ref{lem:Dif}(\ref{itm:Dif-i}), we have $a_k(y)=d_k(x)-d_{k-1}(x)$ for each $k\in\N$, where $d_0(x)\coloneqq 0$, and by \eqref{eq:cylinder-corr}, we have $J_n(y)=\Theta(I_n(x))$. Thus, by Proposition~\ref{prop:Pierce-cylinder} and Proposition~\ref{prop:CF-basic}(\ref{itm:CF-basic-ii}),
\begin{equation}\label{eq:mass}
\begin{aligned}
\mu(J_n(y))&=\len I_n(x)=\frac{1}{Q_n(x)(d_n(x)+1)},\\
\nu(I_n(x))&=\len J_n(y)=\frac{1}{q_n(y)(q_n(y)+q_{n-1}(y))},
\end{aligned}
\end{equation}
where $Q_n(x)\coloneqq d_1(x)\cdots d_n(x)$. Notice that both $(I_n(x))_{n\in\N}$ and $(J_n(y))_{n\in\N}$ are strictly nested sequences of intervals that shrink to $x$ and $y$, respectively.

\begin{proof}[Proof of Theorem~\ref{thm:dimension}]
(\ref{itm:dimension-i}) Let $y\in[0,1]$ and $r>0$. Since $\Theta$ is an increasing homeomorphism, $\Theta^{-1}(B[y,r])$ is an interval with endpoints $\Theta^{-1}(u)$ and $\Theta^{-1}(v)$ for some $u,v\in[0,1]$ with $|u-v|\le2r$. Hence, by Theorem~\ref{thm:quadratic}, we have
\[
\mu(B[y,r])=\len\Theta^{-1}(B[y,r])\le2\sqrt{2r}.
\]

(\ref{itm:dimension-ii}) Let $X$ be the set of irrational $x\in[0,1]$ satisfying both (\ref{lim:L1}) and (\ref{lim:L2}). Clearly, $X$ is a Borel set with $\Leb(X)=1$, and so $\mu(\Theta(X))=1$. Fix $x\in X$ and put $y\coloneqq \Theta(x)$.

By \eqref{eq:mass}, we have $-\log\mu(J_n(y))=\sum_{k=1}^n\log d_k(x)+\log(d_n(x)+1)$. Let $\varepsilon>0$. By (\ref{lim:L1}), there exists $K\in\N$ such that $(1-\varepsilon)k\le\log d_k(x)\le(1+\varepsilon)k$ for all $k\ge K$. Summing over $K\le k\le n$ and dividing by $n^2/2$, we see that both the limit inferior and the limit superior of $\bigl(\sum_{k=1}^n\log d_k(x)\bigr)/(n^2/2)$ lie in $[1-\varepsilon,1+\varepsilon]$. Since $\varepsilon>0$ is arbitrary and $(\log(d_n(x)+1))/n^2\to0$, it follows that
\[
\lim_{n\to\infty}\frac{-\log\mu(J_n(y))}{n^2}=\frac12.
\]
Similarly, Proposition~\ref{prop:CF-basic}(\ref{itm:CF-basic-iii}) gives
\[
\sum_{k=1}^n\log a_k(y)\le\log q_n(y)\le\sum_{k=1}^n\log a_k(y)+n\log2,
\]
and $(\log a_k(y))/k\to1$ as $k\to\infty$ by (\ref{lim:L2}). Hence, by the same argument, $(\log q_n(y))/n^2\to1/2$ as $n\to\infty$. Since $1/(2q_n(y)^2)\le\len J_n(y)\le1/q_n(y)^2$, it follows that
\[
\lim_{n\to\infty}\frac{-\log\len J_n(y)}{n^2}=1.
\]

Let $r>0$ be sufficiently small. Then we can choose $n\in\N$ such that $\len J_{n+1}(y)\le r<\len J_n(y)$. Since $y\in J_{n+1}(y)$ and $\diam J_{n+1}(y)\le r$, we have $J_{n+1}(y)\subseteq B[y,r]$, and so $\mu(B[y,r])\ge\mu(J_{n+1}(y))$. Therefore,
\[
\frac{\log\mu(B[y,r])}{\log r}\le\frac{-\log\mu(J_{n+1}(y))}{-\log\len J_n(y)}=\frac{-\log\mu(J_{n+1}(y))}{(n+1)^2}\cdot\frac{(n+1)^2}{n^2}\cdot\frac{n^2}{-\log\len J_n(y)}.
\]
Since $n\to\infty$ as $r\to0^+$, the right-hand side tends to $(1/2)\cdot1\cdot1=1/2$ as $r\to0^+$. On the other hand, by part (\ref{itm:dimension-i}),
\[
\frac{\log\mu(B[y,r])}{\log r}\ge\frac{\log(2\sqrt2)+(1/2)\log r}{\log r}\to\frac12
\]
as $r\to0^+$. Thus, $\dimloc\mu(y)=1/2$ for every $y\in\Theta(X)$. Since $\mu(\Theta(X))=1$, the assertion follows.

(\ref{itm:dimension-iii}) Let $Y$ be the Borel set of irrational $y\in[0,1]$ satisfying (\ref{lim:L3}), and put $X'\coloneqq \Theta^{-1}(Y)$. Then $\nu(X')=\Leb(Y)=1$. Fix $x\in X'$ and put $y\coloneqq \Theta(x)$. By \eqref{eq:mass}, we have $-\log\nu(I_n(x))=-\log\len J_n(y)\le\log(2q_n(y)^2)$. Since $(\log q_n(y))/n$ converges by (\ref{lim:L3}), it is bounded, and so there exists a constant $C>0$ such that
\[
-\log\nu(I_n(x))\le Cn
\]
for all $n\in\N$. On the other hand, since $d_k(x)\ge k$ for each $k\in\N$, we have $Q_n(x)\ge n!$, and so $\len I_n(x)\le1/(n+1)!$. Let $r>0$ be sufficiently small, and choose $n\in\N$ such that $\len I_{n+1}(x)\le r<\len I_n(x)$. Then $I_{n+1}(x)\subseteq B[x,r]$, and it follows that
\[
0\le\frac{\log\nu(B[x,r])}{\log r}\le\frac{-\log\nu(I_{n+1}(x))}{-\log\len I_n(x)}\le\frac{C(n+1)}{\log((n+1)!)},
\]
which tends to $0$ as $r\to0^+$. Thus $\dimloc\nu(x)=0$.

Finally, the last assertion follows from Proposition~\ref{prop:local-dim}.
\end{proof}

\begin{proof}[Proof of Corollary~\ref{cor:consequences}]
(\ref{itm:consequences-i}) By Theorem~\ref{thm:dimension}, we have $\dimH\mu=1/2<1$ and $\dimH\nu=0<1$. So, by the definition of the Hausdorff dimension of a measure, each of $\mu$ and $\nu$ is carried by a Borel set of Hausdorff dimension strictly less than $1$. Such a set is of null Lebesgue measure, and so $\mu$ and $\nu$ are singular with respect to $\Leb$. Since $\nu([0,x])=\Theta(x)$ and $\mu([0,y])=\Theta^{-1}(y)$, and since the distribution function of a finite measure which is singular with respect to $\Leb$ has derivative $0$ at $\Leb$-almost every point (see, e.g., \cite[Theorem~7.13]{Rud87}), the assertion follows.

(\ref{itm:consequences-ii}) Let $Y_0$ be the set of points $y\in[0,1]$ at which $(\Theta^{-1})'(y)=0$. Then $\Leb(Y_0)=1$ by part (\ref{itm:consequences-i}). Let $y\in Y_0$ and $x\coloneqq \Theta^{-1}(y)$. For $x'\in[0,1]$ with $x'\neq x$, put $y'\coloneqq \Theta(x')$. Then $y'\to y$ as $x'\to x$, and so
\[
\frac{\Theta(x')-\Theta(x)}{x'-x}=\frac{y'-y}{\Theta^{-1}(y')-\Theta^{-1}(y)}\to+\infty
\]
as $x'\to x$. Hence, $\Theta'=+\infty$ on $\Theta^{-1}(Y_0)$, and $\nu(\Theta^{-1}(Y_0))=\Leb(Y_0)=1$.

(\ref{itm:consequences-iii}) The first assertion is a consequence of Theorem~\ref{thm:quadratic} and Proposition~\ref{prop:holder-dim}. Suppose, to the contrary, that $\Theta^{-1}$ is H\"older continuous of exponent $\alpha$ for some $\alpha>1/2$. Then the same argument as in the proof of Theorem~\ref{thm:dimension}(\ref{itm:dimension-i}) gives a constant $C>0$ such that $\mu(B[y,r])\le Cr^\alpha$ for all $y\in[0,1]$ and $r>0$. This implies $\dimloc\mu(y)\ge\alpha>1/2$ whenever the local dimension exists, which contradicts Theorem~\ref{thm:dimension}(\ref{itm:dimension-ii}).

(\ref{itm:consequences-iv}) Suppose, to the contrary, that there are constants $C>0$ and $\alpha>0$ such that $|\Theta(x)-\Theta(x')|\le C|x-x'|^\alpha$ for all $x,x'\in[0,1]$. Then $\nu(B[x,r])=\Leb(\Theta(B[x,r]))\le C(2r)^\alpha$ for all $x\in[0,1]$ and $r>0$. This implies $\dimloc\nu(x)\ge\alpha>0$ whenever the local dimension exists, which contradicts Theorem~\ref{thm:dimension}(\ref{itm:dimension-iii}).
\end{proof}

\begin{remark}\label{rem:nu-carrier}
As mentioned in Remarks~\ref{rems:dimension}(2), we can find a natural set carrying $\nu$. For irrational $x\in[0,1]$ and $y\coloneqq \Theta(x)$, we have $d_n(x)=\sum_{k=1}^na_k(y)\ge n$ for each $n\in\N$. Let $\varepsilon>0$. By Diamond and Vaaler \cite[Corollary~2]{DV86} with $\varphi(n)=n^{1+\varepsilon}$, there exists a Borel set $Y_\varepsilon\subseteq[0,1]$ with $\Leb(Y_\varepsilon)=1$ such that, for each $y\in Y_\varepsilon$,
\[
\sum_{k=1}^na_k(y)\le\frac{2}{\log2}\,n\log n+n^{1+\varepsilon}
\]
for all sufficiently large $n$. Since $(2/\log2)\,n\log n\le n^{1+\varepsilon}$ for all sufficiently large $n$, it follows that, for each irrational $y\in Y_\varepsilon$ and $x\coloneqq\Theta^{-1}(y)$, we have $d_n(x)\le2n^{1+\varepsilon}$ for all sufficiently large $n$, and hence
\[
\limsup_{n\to\infty}\frac{\log d_n(x)}{\log n}\le1+\varepsilon.
\]
Now, put $Y\coloneqq\bigcap_{m\in\N}Y_{1/m}$. Then $\Leb(Y)=1$, since $Y$ is a countable intersection of sets of full measure. Let $y\in Y$ be irrational, and put $x\coloneqq\Theta^{-1}(y)$. Then the above inequality holds with $\varepsilon=1/m$ for every $m\in\N$, and so the limit superior is at most $1$. On the other hand, the limit inferior is at least $1$ because $d_n(x)\ge n$ for each $n\in\N$. Therefore,
\[
\lim_{n\to\infty}\frac{\log d_n(x)}{\log n}=1
\]
for every $x\in\Theta^{-1}(Y\setminus\Q)$. Notice that $\Theta^{-1}(Y\setminus\Q)$ consists of irrational points, since $\Theta$ maps rationals onto rationals, and that $\nu(\Theta^{-1}(Y\setminus\Q))=\Leb(Y\setminus\Q)=1$ by the definition of $\nu$. Thus, $\nu$ is carried by the set
\[
E\coloneqq \left\{x\in(0,1]\setminus\Q:\lim_{n\to\infty}\frac{\log d_n(x)}{\log n}=1\right\},
\]
which is of null Hausdorff dimension by \cite[Theorem~1.1]{Ahn24} (with $\phi(n)=\log(n+1)$ there). This gives another proof of $\dimH\nu=0$.
\end{remark}

\subsection{Differentiability}

\begin{proof}[Proof of Theorem~\ref{thm:rational}]
Let $r=\langle d_1,\dots,d_n\rangle_P\in(0,1)$, so that $\sigma\coloneqq (d_1,\dots,d_n)\in\Sigma_{\mathrm P,n}\cap\SP^{\mathrm{re}}$. Write $\Theta(r)=[a_1,\dots,a_n]$ with $a_k\coloneqq d_k-d_{k-1}$, where $d_0\coloneqq 0$, and let $p_k/q_k$, $0\le k\le n$, be its convergents. Put $s\coloneqq (-1)^n$. We will show that the one-sided derivative of $\Theta$ at $r$ from the side $s$ (i.e., from the right if $s=1$ and from the left if $s=-1$) is $Q_n/q_n^2$, and the one from the side $-s$ is $Q_n(d_n-1)/q_n^2$.

{\sc The side $s$.} By Proposition~\ref{prop:Pierce-cylinder}, $I_\sigma=\{r+st/Q_n:0\le t<1/(d_n+1)\}$, which is a one-sided neighborhood of $r$ on the side $s$. Let $x\in I_\sigma\setminus\{r\}$. Then $t=T^n(x)>0$, and the Pierce digits of $t$ are $d_{n+1}(x),d_{n+2}(x),\dots$, all of which are larger than $d_n$. By Lemma~\ref{lem:Dif}(\ref{itm:Dif-i}), the continued fraction digits of $\Theta(x)$ after the $n$th one are $d_{n+1}(x)-d_n,\,d_{n+2}(x)-d_{n+1}(x),\dots$. This means that $u\coloneqq S^n(\Theta(x))=F_{d_n}(t)$. So, by Proposition~\ref{prop:CF-basic}(\ref{itm:CF-basic-ii}), $\Theta(x)-\Theta(r)=su/(q_n(q_n+q_{n-1}u))$. Therefore,
\[
\frac{\Theta(x)-\Theta(r)}{x-r}=\frac{Q_n}{q_n^2}\cdot\frac{F_{d_n}(t)}{t}\cdot\frac{q_n}{q_n+q_{n-1}F_{d_n}(t)}\to\frac{Q_n}{q_n^2}
\]
as $x\to r$ from the side $s$, because $t\to0$, $F_{d_n}(t)\to0$, and $F_{d_n}(t)/t\to1$ by Lemma~\ref{lem:ratio}(\ref{itm:ratio-iii}).

{\sc The side $-s$.} Let $\sigma'\coloneqq (d_1,\dots,d_{n-1},d_n-1,d_n)$. If $n\ge2$, then $d_n-1\ge d_{n-1}+1$; if $n=1$, then $d_1\ge2$ because $r<1$. So, $\sigma'\in\Sigma_{\mathrm P,n+1}$, $\sigma'$ is not realizable, and $\phiP(\sigma')=r$ by \eqref{eq:Pierce-two-reps}. Put $Q'\coloneqq Q(\sigma')=Q_n\cdot(d_n-1)$. By Proposition~\ref{prop:Pierce-cylinder}, $I_{\sigma'}=\{r-st/Q':0<t<1/(d_n+1)\}$, which is a punctured one-sided neighborhood of $r$ on the side $-s$. By definition, $\Dif_0(\sigma')=(a_1,\dots,a_{n-1},a_n-1,1)$, where $a_n-1=(d_n-1)-d_{n-1}\in\N$. Its convergent denominators $q'_k$ satisfy $q'_k=q_k$ for $k\le n-1$, $q'_n=q_n-q_{n-1}$, and $q'_{n+1}=q'_n+q_{n-1}=q_n$, and its value is $[a_1,\dots,a_{n-1},a_n-1,1]=[a_1,\dots,a_n]=\Theta(r)$. As before, for $x\in I_{\sigma'}$, we have $S^{n+1}(\Theta(x))=F_{d_n}(t)=:u$, and so
\[
\Theta(x)-\Theta(r)=\frac{-su}{q_n(q_n+q'_nu)}
\]
and
\[
\frac{\Theta(x)-\Theta(r)}{x-r}=\frac{Q'}{q_n^2}\cdot\frac{F_{d_n}(t)}{t}\cdot\frac{q_n}{q_n+q'_nF_{d_n}(t)}\to\frac{Q_n(d_n-1)}{q_n^2}
\]
as $x\to r$ from the side $-s$.

{\sc The endpoints.} We have $\Theta(x)/x\to1$ as $x\to0^+$ by Lemma~\ref{lem:ratio}(\ref{itm:ratio-iii}). For $x\in(1/2,1)$, Lemma~\ref{lem:recursion} with $D=0$ and $j=1$ gives $\Theta(x)=1/(1+F_1(1-x))$. Hence,
\[
\frac{1-\Theta(x)}{1-x}=\frac{F_1(1-x)}{(1-x)(1+F_1(1-x))}\to1
\]
as $x\to1^-$, again by Lemma~\ref{lem:ratio}(\ref{itm:ratio-iii}). This completes the proof.
\end{proof}

\begin{proof}[Proof of Theorem~\ref{thm:irrational}]
Write $d_n\coloneqq d_n(x)$, $Q_n\coloneqq Q_n(x)$, $y\coloneqq \Theta(x)$, $a_n\coloneqq a_n(y)=d_n-d_{n-1}$, and $q_n\coloneqq q_n(y)$. Suppose, to the contrary, that $\Theta'(x)=L$ for some $L\in(0,\infty)$.

{\sc Step 1.} Since $x$ is irrational, $x$ lies in the interior of $I_n(x)$, whose endpoints are rational, for each $n\in\N$. Let $x_n^-<x<x_n^+$ be the endpoints of the closure of $I_n(x)$. Then the slope $(\Theta(x_n^+)-\Theta(x_n^-))/(x_n^+-x_n^-)$ is a convex combination of the two difference quotients of $\Theta$ over $[x_n^-,x]$ and $[x,x_n^+]$, and so it tends to $L$ as $n\to\infty$. By \eqref{eq:cylinder-corr}, this slope equals
\[
\rho_n\coloneqq \frac{\len J_n(y)}{\len I_n(x)}=\frac{Q_n(d_n+1)}{q_n(q_n+q_{n-1})}.
\]
Hence $\rho_n\to L$, and so $\rho_n/\rho_{n-1}\to1$ as $n\to\infty$.

{\sc Step 2.} We claim that $a_n\to\infty$ as $n\to\infty$. Using $Q_n=d_nQ_{n-1}$ and $q_n=a_nq_{n-1}+q_{n-2}$, and putting $\eta_n\coloneqq q_{n-2}/q_{n-1}\in[0,1]$, we find that
\[
\frac{\rho_n}{\rho_{n-1}}=\frac{d_n(d_n+1)}{d_{n-1}+1}\cdot\frac{1+\eta_n}{(a_n+\eta_n)(a_n+\eta_n+1)}
\]
for each $n\ge2$. If $a_n\le M$ for infinitely many $n$, for some $M\in\N$, then the right-hand side is at least $d_n/((M+1)(M+2))$ for such $n$, which tends to $\infty$ since $d_n\to\infty$. This contradicts $\rho_n/\rho_{n-1}\to1$, and the claim is proved.

{\sc Step 3.} Put $r_n\coloneqq \phiP(d_1,\dots,d_n)$ for each $n\in\N$. Then $r_n\neq x$, $r_n\to x$ as $n\to\infty$, and $\Theta(r_n)=[a_1,\dots,a_n]=p_n/q_n$ by Lemma~\ref{lem:Dif}(\ref{itm:Dif-ii}). By \eqref{eq:Pierce-tail} and Proposition~\ref{prop:CF-basic}(\ref{itm:CF-basic-ii}), we have
\[
x-r_n=\frac{(-1)^nt_n}{Q_n}\qquad\text{and}\qquad\Theta(x)-\Theta(r_n)=\frac{(-1)^nu_n}{q_n(q_n+q_{n-1}u_n)},
\]
where $t_n\coloneqq T^n(x)\in(1/(d_{n+1}+1),1/d_{n+1})$ and $u_n\coloneqq S^n(y)\in(1/(a_{n+1}+1),1/a_{n+1})$; here, the intervals are open since $x$ and $y$ are irrational. Put $R_n\coloneqq (\Theta(x)-\Theta(r_n))/(x-r_n)$. Then, using $d_{n+1}=d_n+a_{n+1}$, we obtain
\[
0<\frac{R_n}{\rho_n}=\frac{u_n}{t_n}\cdot\frac{q_n+q_{n-1}}{(d_n+1)(q_n+q_{n-1}u_n)}\le\frac{d_{n+1}+1}{a_{n+1}}\cdot\frac{2}{d_n+1}=2\left(\frac{1}{a_{n+1}}+\frac{1}{d_n+1}\right).
\]
By Step 2, the right-hand side tends to $0$ as $n\to\infty$, and so $R_n=\rho_n\cdot(R_n/\rho_n)\to L\cdot0=0$. But $R_n\to\Theta'(x)=L>0$ since $r_n\to x$, which is a contradiction. This completes the proof.
\end{proof}

\begin{proof}[Proof of Corollary~\ref{cor:classification}]
By Theorem~\ref{thm:rational}, the two one-sided derivatives of $\Theta$ at a rational $r=\langle d_1,\dots,d_n\rangle_P\in(0,1)$ differ by the factor $d_n-1$, so they agree if and only if $d_n=2$. If $n\ge2$, then $d_n\ge d_{n-1}+2\ge3$ by Proposition~\ref{prop:Pierce-basic}(\ref{itm:Pierce-basic-ii}). Thus $d_n=2$ forces $n=1$, i.e., $r=1/2$, in which case $Q_1=q_1=2$ and both one-sided derivatives are equal to $1/2$. Combining this with Theorem~\ref{thm:irrational}, the corollary follows.
\end{proof}

\subsection{Dimension of the set of infinite derivatives}

\begin{proof}[Proof of Theorem~\ref{thm:infinite-dimension}]
For irrational $x\in(0,1)$, write $y\coloneqq\Theta(x)$, $d_n\coloneqq d_n(x)$, $a_n\coloneqq a_n(y)=d_n-d_{n-1}$, $Q_n\coloneqq d_1\cdots d_n$, and $q_n\coloneqq q_n(y)$, where $d_0\coloneqq0$. As in the proof of Theorem~\ref{thm:irrational}, put
\[
\rho_n(x)\coloneqq\frac{\len J_n(y)}{\len I_n(x)}
=\frac{Q_n(d_n+1)}{q_n(q_n+q_{n-1})}.
\]
We establish the upper and lower bounds separately.

{\sc Upper bound.} By Theorem~\ref{thm:rational}, every point of $E_\infty$ is irrational. If $x\in E_\infty$, then $\rho_n(x)\to+\infty$ as $n\to\infty$. Indeed, the slope over $I_n(x)$ is a convex combination of the two difference quotients from $x$ to the endpoints of this interval, and both quotients tend to $+\infty$. This is the same argument as in Step~1 of the proof of Theorem~\ref{thm:irrational}.

For each $n\in\N$, let
\[
\mathcal W_n\coloneqq
\{\sigma\in\Sigma_{\mathrm P,n}:
\len I_\sigma\le\len J_{\Dif_0(\sigma)}\}.
\]
It follows that, for each $N\in\N$, the collection of intervals $I_\sigma$ with $\sigma\in\mathcal W_n$ and $n\ge N$ covers $E_\infty$. By Proposition~\ref{prop:Pierce-cylinder}, every Pierce fundamental interval of order $n$ has length at most
\[
\kappa_n\coloneqq\frac1{(n+1)!},
\]
since its digits satisfy $\sigma_k\ge k$ for $1\le k\le n$.

Fix $s>1/2$, and choose $t$ with $1/2<t<s$. By Proposition~\ref{prop:CF-basic}(\ref{itm:CF-basic-iii}), we have $q_n\ge a_1\cdots a_n$. Since $\Dif_0$ maps $\Sigma_{\mathrm P,n}$ bijectively onto $\Sigma_{\mathrm C,n}$, which we identify with $\N^n$, we obtain
\[
\begin{split}
\sum_{\sigma\in\mathcal W_n}(\len I_\sigma)^s
&\le\kappa_n^{s-t}
\sum_{\sigma\in\mathcal W_n}(\len J_{\Dif_0(\sigma)})^t\\
&\le\kappa_n^{s-t}
\sum_{(a_1,\dots,a_n)\in\N^n}(a_1\cdots a_n)^{-2t}
=\kappa_n^{s-t}Z_t^n,
\end{split}
\]
where $Z_t\coloneqq\sum_{j\in\N}j^{-2t}<\infty$. The series
\[
\sum_{n\in\N}\kappa_n^{s-t}Z_t^n
=\sum_{n\in\N}\frac{Z_t^n}{((n+1)!)^{s-t}}
\]
converges by the ratio test, since $s-t>0$ and
\[
\frac{Z_t^{n+1}/((n+2)!)^{s-t}}
     {Z_t^n/((n+1)!)^{s-t}}
=\frac{Z_t}{(n+2)^{s-t}}\to0
\]
as $n\to\infty$. For each $N\in\N$, the above covering satisfies
\[
\sup\{\diam I_\sigma:n\ge N,\,\sigma\in\mathcal W_n\}
\le\kappa_N\to0
\qquad\text{as }N\to\infty,
\]
and
\[
\sum_{n=N}^{\infty}
\sum_{\sigma\in\mathcal W_n}(\diam I_\sigma)^s
\le
\sum_{n=N}^{\infty}
\frac{Z_t^n}{((n+1)!)^{s-t}}
\to0
\qquad\text{as }N\to\infty.
\]
Hence, the $s$-dimensional Hausdorff measure of $E_\infty$ equals $0$, and so
$\dimH E_\infty\le s$. Since $s>1/2$ was arbitrary, this proves the upper bound.

{\sc Lower bound.} Fix $0<\alpha<1$, and put $u_n\coloneqq2n^\alpha$ for each $n\in\N$. Consider the set
\[
K_\alpha\coloneqq
\{x\in(0,1)\setminus\Q:
nu_n<d_n(x)\le(n+1)u_n\text{ for all }n\in\N\}.
\]
The sequence $(u_n)_{n\in\N}$ is non-decreasing and satisfies $u_n\ge2$. By \cite[Lemma~3.2]{Ahn24}, we have
\[
\dimH K_\alpha=\frac{\alpha}{\alpha+1}.
\]
In fact, the parameter $\eta$ in that lemma equals $1/\alpha$, since
\[
\frac{n\log n+\log u_{n+1}}{\sum_{k=1}^n\log u_k}
=\frac{n\log n+\log2+\alpha\log(n+1)}
{n\log2+\alpha\log(n!)}
\to\frac1\alpha
\]
as $n\to\infty$.

We show that $K_\alpha\subseteq E_\infty$. Fix $x\in K_\alpha$, and use the notation from the beginning of the proof. The defining inequalities give $2n^{\alpha+1}<d_n\le2(n+1)n^\alpha\le4n^{\alpha+1}$ for each $n\in\N$. Let $n\ge2$. By the mean value theorem, since $0<\alpha<1$, we have
\[
a_n=d_n-d_{n-1}>2n\bigl(n^\alpha-(n-1)^\alpha\bigr)\ge2n\cdot\alpha n^{\alpha-1}=2\alpha n^\alpha
\]
and
\[
a_n\le2(n+1)n^\alpha-2(n-1)^{\alpha+1}=2\bigl(n^{\alpha+1}-(n-1)^{\alpha+1}\bigr)+2n^\alpha\le2(\alpha+2)n^\alpha\le6n^\alpha.
\]
Also, $2<a_1=d_1\le4$. Hence, $\log Q_n>(\alpha+1)\log(n!)+n\log2$, and, by Proposition~\ref{prop:CF-basic}(\ref{itm:CF-basic-iii}),
\[
\log q_n\le\sum_{k=1}^n\log a_k+n\log2\le\alpha\log(n!)+n\log12.
\]
Since $\rho_n(x)\ge Q_n/(2q_n^2)$, we obtain
\begin{equation}\label{eq:infinite-slope}
\log\rho_n(x)\ge(1-\alpha)\log(n!)-2n\log12
\end{equation}
for each $n\in\N$.

To pass from the slopes over fundamental intervals to the derivative, we need an estimate for the distance of $y$ from their endpoints. For each $m\in\N$, let $L_m\coloneqq\overline{J_m(y)}$ and $\delta_m\coloneqq\operatorname{dist}(y,\partial L_m)$. If $a_{m+1}\ge2$, put $w\coloneqq S^m(y)$, $a\coloneqq a_{m+1}$, and $r\coloneqq q_{m-1}/q_m$. Then $0\le r\le1$ and $1/(a+1)<w<1/a$. By Proposition~\ref{prop:CF-basic}(\ref{itm:CF-basic-ii}), the distances of $y$ from the two endpoints of $L_m$, divided by $\len J_m(y)$, are
\[
\frac{w(1+r)}{1+rw}
\quad\text{and}\quad
\frac{1-w}{1+rw}.
\]
The first is at least $w>1/(a+1)$, while the second is at least $(1-1/a)/(1+1/a)=(a-1)/(a+1)\ge1/(a+1)$. Thus,
\[
\delta_m\ge\frac{\len J_m(y)}{a_{m+1}+1}
\qquad\text{if }a_{m+1}\ge2.
\]

Put $G\coloneqq\Theta^{-1}$ and $L_0\coloneqq[0,1]$. For $z\in[0,1]\setminus\{y\}$ sufficiently close to $y$, let $n$ be the largest integer such that $z\in L_n$. Such an integer exists because the intervals $L_n$ decrease to $\{y\}$. Also, $n\to\infty$ as $z\to y$, since $y$ lies in the interior of each $L_n$. We have $|z-y|\ge\delta_{n+1}$, while $|G(z)-G(y)|\le\len I_n(x)$ by \eqref{eq:cylinder-corr} and continuity of $G$. Since $a_m>2\alpha m^\alpha$ for all $m\ge2$,
we have $a_m\to\infty$ as $m\to\infty$. Hence the preceding estimate for $\delta_{n+1}$ applies for all sufficiently large $n$. Moreover,
\[
\frac{\len J_n(y)}{\len J_{n+1}(y)}
=\frac{q_{n+1}(q_{n+1}+q_n)}{q_n(q_n+q_{n-1})}
\le(a_{n+1}+1)(a_{n+1}+2).
\]
Therefore,
\[
\begin{split}
\frac{|G(z)-G(y)|}{|z-y|}
&\le\frac{\len I_n(x)}{\delta_{n+1}}
\le(a_{n+2}+1)\frac{\len I_n(x)}{\len J_{n+1}(y)}\\
&\le
\frac{(a_{n+2}+1)(a_{n+1}+1)(a_{n+1}+2)}
{\rho_n(x)}.
\end{split}
\]
Since $6t+2\le8t$ for $t\ge1$, the numerator is at most $8^3(n+2)^{3\alpha}$, whereas $\rho_n(x)\ge(n!)^{1-\alpha}/144^n$ by \eqref{eq:infinite-slope}. Since $(n!)^{1-\alpha}/(144^n(n+2)^{3\alpha})\to\infty$ as $n\to\infty$, the right-hand side tends to $0$, and so $G'(y)=0$. Taking reciprocal difference quotients, as in the proof of Corollary~\ref{cor:consequences}(\ref{itm:consequences-ii}), we obtain $\Theta'(x)=+\infty$. This proves $K_\alpha\subseteq E_\infty$. Consequently,
\[
\dimH E_\infty\ge\dimH K_\alpha=\frac{\alpha}{\alpha+1}.
\]
Letting $\alpha\to1^-$ yields the lower bound $1/2$, and completes the proof.
\end{proof}

\section{Concluding remark}\label{sec:remarks}

Our map $\Theta$ compares Pierce expansions with continued fraction expansions by taking the gaps between consecutive Pierce digits as continued fraction digits. We conclude this paper with a related comparison of the two expansions of the same irrational number $x$, namely, how many continued fraction digits of $x$ are determined by the first $n$ Pierce digits of $x$.

For irrational $x\in[0,1]$, put $J_0(x)\coloneqq [0,1]$ and define
\[
m_n(x)\coloneqq \max\{m\in\N_0:I_n(x)\subseteq J_m(x)\}
\]
for each $n\in\N$. The intervals used in \cite{BCRS23} are open. To match their convention, put
\[
\widetilde m_n(x)\coloneqq
\max\{m\in\N_0:I_n(x)^\circ\subseteq J_m(x)^\circ\},
\]
where $A^\circ$ denotes the interior of $A$. Then
\[
0\le\widetilde m_n(x)-m_n(x)\le2.
\]
Indeed, $m_n(x)\le\widetilde m_n(x)$ is immediate. For $m\ge2$, the continued fraction cylinder formulas in Proposition~\ref{prop:CF-basic}(\ref{itm:CF-basic-ii}) give $\overline{J_m(x)}\subseteq J_{m-2}(x)$. Hence, if $\widetilde m_n(x)\ge2$, then
\[
I_n(x)\subseteq\overline{I_n(x)^\circ}
\subseteq\overline{J_{\widetilde m_n(x)}(x)}
\subseteq J_{\widetilde m_n(x)-2}(x),
\]
which gives the other inequality. If $\widetilde m_n(x)\le1$, it follows from $m_n(x)\ge0$.

By (\ref{lim:L1}) and (\ref{lim:L3}), we have
\[
\lim_{n\to\infty}\frac{-\log\len I_n(x)}{n^2/2}=1
\qquad\text{and}\qquad
\lim_{m\to\infty}\frac{-\log\len J_m(x)}{\pi^2m/(6\log2)}=1
\]
for $\Leb$-almost every $x$. This means that the Pierce and continued fraction partitions are almost everywhere log-balanced, in the sense of Berth\'e, Cesaratto, Rotondo, and Safe \cite{BCRS23}, with weights
\[
f_1(n)\coloneqq \frac{n^2}{2}\qquad\text{and}\qquad f_2(m)\coloneqq \frac{\pi^2m}{6\log2},
\]
respectively. Clearly, $f_2$ is non-decreasing with $f_2(m+1)/f_2(m)\to1$ as $m\to\infty$. Moreover, we have $\sum_{n\in\N}e^{-\delta f_1(n)}<\infty$ for every $\delta>0$. Hence, \cite[Theorem~6.8]{BCRS23}, applied to the interiors of the fundamental intervals, gives
\[
\lim_{n\to\infty}\frac{f_2(\widetilde m_n(x))}{f_1(n)}=1
\]
for $\Leb$-almost every $x$. Since $0\le\widetilde m_n(x)-m_n(x)\le2$, it follows that
\[
\lim_{n\to\infty}\frac{m_n(x)}{n^2}=\frac{3\log2}{\pi^2}
\]
for $\Leb$-almost every $x$.




\begin{thebibliography}{99}

\bibitem{Ahn23} M.~W.~Ahn, {\it On the error-sum function of Pierce expansions}, J.~Fractal Geom. {\bf 10} (2023), 389--421.

\bibitem{Ahn24} M.~W.~Ahn, {\it Hausdorff dimensions in Pierce expansions}, Acta Arith. {\bf 215} (2024), 115--160.

\bibitem{Ahn25} M.~W.~Ahn, {\it Continuity of the continued fraction mapping revisited}, Arch.~Math. {\bf 124} (2025), 395--405.

\bibitem{ABPT07} S.~Albeverio, O.~Baranovskyi, M.~Pratsiovytyi, \& G.~Torbin, {\it The Ostrogradsky series and related Cantor-like sets}, Acta Arith. {\bf 130} (2007), 215--230.

\bibitem{AT15} S.~Albeverio \& G.~Torbin, {\it Ostrogradsky--Sierpi\'nski--Pierce expansion: dynamical systems, probability theory and fractal geometry points of view}, preprint (2015), \href{https://arxiv.org/abs/1506.04355}{arXiv:1506.04355}.

\bibitem{BCRS23} V.~Berth\'e, E.~Cesaratto, P.~Rotondo, \& M.~D.~Safe, {\it Lochs-type theorems beyond positive entropy}, Monatsh. Math. {\bf 200} (2023), 737--779.

\bibitem{dADF22} E.~de Amo, M.~D\'iaz Carrillo, \& J.~Fern\'andez-S\'anchez, {\it A study of topological conjugacies between alternate representation systems}, Rev. R. Acad. Cienc. Exactas F\'is. Nat. Ser. A Mat. RACSAM {\bf 116} (2022), article no. 136.

\bibitem{Den38} A.~Denjoy, {\it Sur une fonction r\'eelle de Minkowski}, J.~Math. Pures Appl. (9) {\bf 17} (1938), 105--151.

\bibitem{DV86} H.~G.~Diamond \& J.~D.~Vaaler, {\it Estimates for partial sums of continued fraction partial quotients}, Pacific J.~Math. {\bf 122} (1986), 73--82.

\bibitem{Fal97} K.~J.~Falconer, Techniques in Fractal Geometry, John Wiley \& Sons, Chichester, 1997.

\bibitem{Fal14} K.~J.~Falconer, Fractal Geometry: Mathematical Foundations and Applications, 3rd edn., John Wiley \& Sons, Chichester, 2014.

\bibitem{FVPD12} J.~Fern\'andez-S\'anchez, P.~Viader, J.~Parad\'is, \& M.~D\'iaz Carrillo, {\it A singular function with a non-zero finite derivative}, Nonlinear Anal. {\bf 75} (2012), 5010--5014.

\bibitem{FVPD14} J.~Fern\'andez-S\'anchez, P.~Viader, J.~Parad\'is, \& M.~D\'iaz Carrillo, {\it A singular function with a non-zero finite derivative on a dense set}, Nonlinear Anal. {\bf 95} (2014), 703--713.

\bibitem{Goo41} I.~J.~Good, {\it The fractional dimensional theory of continued fractions}, Proc. Cambridge Philos. Soc. {\bf 37} (1941), 199--228.

\bibitem{IK02} M.~Iosifescu \& C.~Kraaikamp, Metrical Theory of Continued Fractions, Math. Appl., Vol. 547, Kluwer Academic Publishers, Dordrecht, 2002.

\bibitem{Khi64} A.~Ya.~Khinchin, Continued Fractions, University of Chicago Press, Chicago, 1964.

\bibitem{Lev36} P.~L\'evy, {\it Sur le d\'eveloppement en fraction continue d'un nombre choisi au hasard}, Compositio Math. {\bf 3} (1936), 286--303.

\bibitem{LLL25} C.~Long, Y.~Liang, \& L.~Lu, {\it Hausdorff dimension of the limsup and liminf sets in Pierce expansions}, Fractals {\bf 33} (2025), article no. 2550029.

\bibitem{LLS26} L.~Lu, C.~Long, \& L.~Shang, {\it Limit theorems and fractal properties of digit gaps in Pierce expansions}, J.~Math.~Anal.~Appl. {\bf 555} (2026), article no. 130019.

\bibitem{PVB99} J.~Parad\'is, P.~Viader, \& L.~Bibiloni, {\it A total order in $(0,1]$ defined through a `next' operator}, Order {\bf 16} (1999), 207--220.

\bibitem{PVB01} J.~Parad\'is, P.~Viader, \& L.~Bibiloni, {\it The derivative of Minkowski's ${?}(x)$ function}, J.~Math.~Anal.~Appl. {\bf 253} (2001), 107--125.

\bibitem{Pie29} T.~A.~Pierce, {\it On an algorithm and its use in approximating roots of algebraic equations}, Amer. Math. Monthly {\bf 36} (1929), 523--525.

\bibitem{PB05} M.~V.~Pratsiovytyi \& O.~M.~Baranovskyi, {\it Properties of distributions of random variables with independent differences of consecutive elements of the Ostrogradski\u{\i} series}, Theory Probab. Math. Statist. (2005), no.~70, 147--160.

\bibitem{Rud87} W.~Rudin, Real and Complex Analysis, 3rd edn., McGraw-Hill, New York, 1987.

\bibitem{Sal43} R.~Salem, {\it On some singular monotonic functions which are strictly increasing}, Trans. Amer. Math. Soc. {\bf 53} (1943), 427--439.

\bibitem{Sha86} J.~O.~Shallit, {\it Metric theory of Pierce expansions}, Fibonacci Quart. {\bf 24} (1986), 22--40.

\bibitem{VBP99} P.~Viader, L.~Bibiloni \& J.~Parad\'is, {\it On a problem of Alfr\'ed R\'enyi}, Acta Arith. {\bf 91} (1999), 107--115.

\end{thebibliography}
\end{document}